\newif
\ifconfver
\confvertrue        
\ifconfver
    \documentclass[10pt,journal,final,twoside]{IEEEtran}
\else
    \documentclass[11pt,draftcls,onecolumn]{IEEEtran}
\fi

\usepackage[utf8]{inputenc}
\usepackage{xspace,empheq,fancybox,amsmath,amssymb,amsfonts,graphicx,epstopdf,epsfig,syntonly,times,amsthm} 
\usepackage{psfrag,color,bm,array,tabularx}
\usepackage{cite,url,footnote,xspace,syntonly,algorithmic}
\usepackage{verbatim,multirow}
\usepackage[linesnumbered,ruled,vlined]{algorithm2e}
\usepackage[T1]{fontenc}
\usepackage{textcomp}
\usepackage{xcolor}
\usepackage{mathtools}
\usepackage{makecell}
\usepackage{bbm, booktabs, makecell, subcaption, titlesec}
\newcolumntype{M}[1]{>{\centering\arraybackslash}m{#1}}
\newcolumntype{N}{@{}m{0pt}@{}}
\DeclareMathOperator*{\argmax}{\arg\max}
\def\BibTeX{{\rm B\kern-.05em{\sc i\kern-.025em b}\kern-.08em
    T\kern-.1667em\lower.7ex\hbox{E}\kern-.125emX}}
\newcommand{\KB}{\color{black}{}}    
\newcommand{\KBB}{\color{black}{}}  

\input{mysymbol.sty}
\allowdisplaybreaks

\IEEEoverridecommandlockouts
\title{\LARGE \bf Reinforcement Learning-based Output Structured Feedback for Distributed Multi-Area Power System Frequency Control}
\author{Kyung-Bin Kwon, Sayak Mukherjee, Hao Zhu and Thanh Long Vu
\thanks{K. Kwon and H. Zhu are with the Department of Electrical and Computer
Engineering, The University of Texas at Austin, 2501 Speedway, Austin,
TX, 78712, USA; Emails: \{kwon8908kr, haozhu\}@utexas.edu.}
\thanks{S. Mukherjee and T.L. Vu are with the Optimization and Control Group, Pacific Northwest National Laboratory, Richland, WA 99352, USA, Emails: \{sayak.mukherjee, thanhlong.vu\}@pnnl.gov.}%
}

\begin{document}
\maketitle

\begin{abstract}
 Load frequency control (LFC) is a key factor to maintain the stable frequency in multi-area power systems. As the modern power systems evolve from centralized to distributed paradigm, LFC needs to consider the peer-to-peer (P2P) based scheme that considers limited information from the information-exchange graph for the generator control of each interconnected area. This paper aims to solve a data-driven constrained LQR problem with mean-variance risk constraints and output structured feedback, and applies this framework to  solve the LFC problem in multi-area power systems. By reformulating the constrained optimization problem into a minimax problem, the stochastic gradient descent max-oracle (SGDmax) algorithm with zero-order policy gradient (ZOPG) is adopted to find the optimal feedback gain from the learning, while guaranteeing convergence. In addition, to improve the adaptation of the proposed learning method to new or varying models, we construct an emulator grid that approximates the dynamics of a physical grid and performs training based on this model. Once the feedback gain is obtained from the emulator grid, it is applied to the physical grid with a robustness test to check whether the controller from the approximated emulator applies to the actual system. Numerical tests show that the obtained feedback controller can successfully control the frequency of each area, while mitigating the uncertainty from the loads, with reliable robustness that ensures the adaptability of the obtained feedback gain to the actual physical grid. 
\end{abstract}

\section{Introduction} \label{sec:IN}

Load frequency control (LFC) is one of the most important control problems in power system operations. The objective of LFC is to maintain the frequency of each area in an interconnected power system by adjusting the output of generators with automatic generation control (AGC) regulator or excitation controller \cite{shashi}. The LFC has been studied in various research to cope with conventional generators \cite{wen}, distributed energy resources \cite{zhu}, and electric vehicles (EVs) \cite{falahati}. However, most research considers a centralized framework, implying that a centralized dispatch center controls all generators \cite{kundur}. As a modern power system relies more on distributed generation (DG) and requires resilience to cyberattacks, the existing centralized control paradigm with one-point-failure faces many challenges. Distributed or decentralized frameworks can be considerable frameworks to improve the stability and security of the power system. Consequently, there is a lot of research that considers distributed or decentralized LFC problems, assuming the limited information exchanges among the interconnected areas \cite{ting, zhao, wang, andrade}. Especially in the case of peer-to-peer (P2P) based LFC, generation control in one area is determined by only the information from the areas connected in the information-exchange graph \cite{ting}.

The general LFC problem can be represented as a linear quadratic regulator (LQR) problem, which minimizes the frequency deviations and other factors such as power outputs, power inflow between interconnected areas, and control efforts such as AGC control signal \cite[Ch.~2]{bevrani}. When the model is known, this problem can be easily solved by adopting the Algebraic Riccati equation (ARE) \cite{ku} or by applying gradient-based methods \cite{bu,li}. However, finding the feedback gain will become complicated as we consider the model uncertainties, the constraints on the optimization problem, or the structure of the feedback gain. This paper aims to solve the P2P-based LFC problem, while dealing with these three challenges at once by considering uncertainty from the environment, the risk constraints, and structured output feedback. 



In particular, reinforcement learning (RL) will be applied to solve this problem. There are two advantages that RL has \cite{sutton}. First, RL is model-free learning, i.e., we do not need to consider the model parameters. Instead, we will generate the trajectory, observe the reward and update the controller toward increasing the total reward. Second, RL is data-driven learning, i.e., instead of generating a probability distribution for the uncertainty and performing the Monte-Carlo method, we can directly use the data gathered from the grid and use it to train the controller. Recent research works such as \cite{sayak_arxiv, eth_distributed_learning, mukh_tac, KB} look into various distributed aspects of RL. 

The main contributions of the paper are three-fold. First, we formulate the constrained LQR problem for LFC considering the mean-variance risk constraint and output structured feedback. Second, an RL-based Stochastic Gradient Descent max-oracle (SGDmax) algorithm with a zero-order policy gradient (ZOPG) is proposed to solve the formulated problem and find the feedback controller. Here, the learning process is applied to the emulator grid which approximates the dynamics of the physical grid. Last, the robustness test on the obtained feedback controller is checked to ensure that the controller trained from the emulator grid has enough robustness that can be directly applied to the physical grid. 

The rest of this paper is organized as follows. In Section~\ref{sec:PS}, problem setup with two models is introduced along with the mathematical formulation of the constrained LQR problem including the objective function, risk constraints, and linear dynamics. Section~\ref{sec:ALG} introduces the SGDmax algorithm with ZOPG to estimate the gradient and find a feedback gain by a gradient-descent method. Section~\ref{sec:SR} shows the effectiveness of the proposed algorithm by adopting it to a load frequency control (LFC) problem. Finally, the paper is wrapped up in Section~\ref{sec:CON}. 

\section{Problem Description} \label{sec:PS}
\begin{figure}[t]
	\centering
	\includegraphics[width=0.9\linewidth]{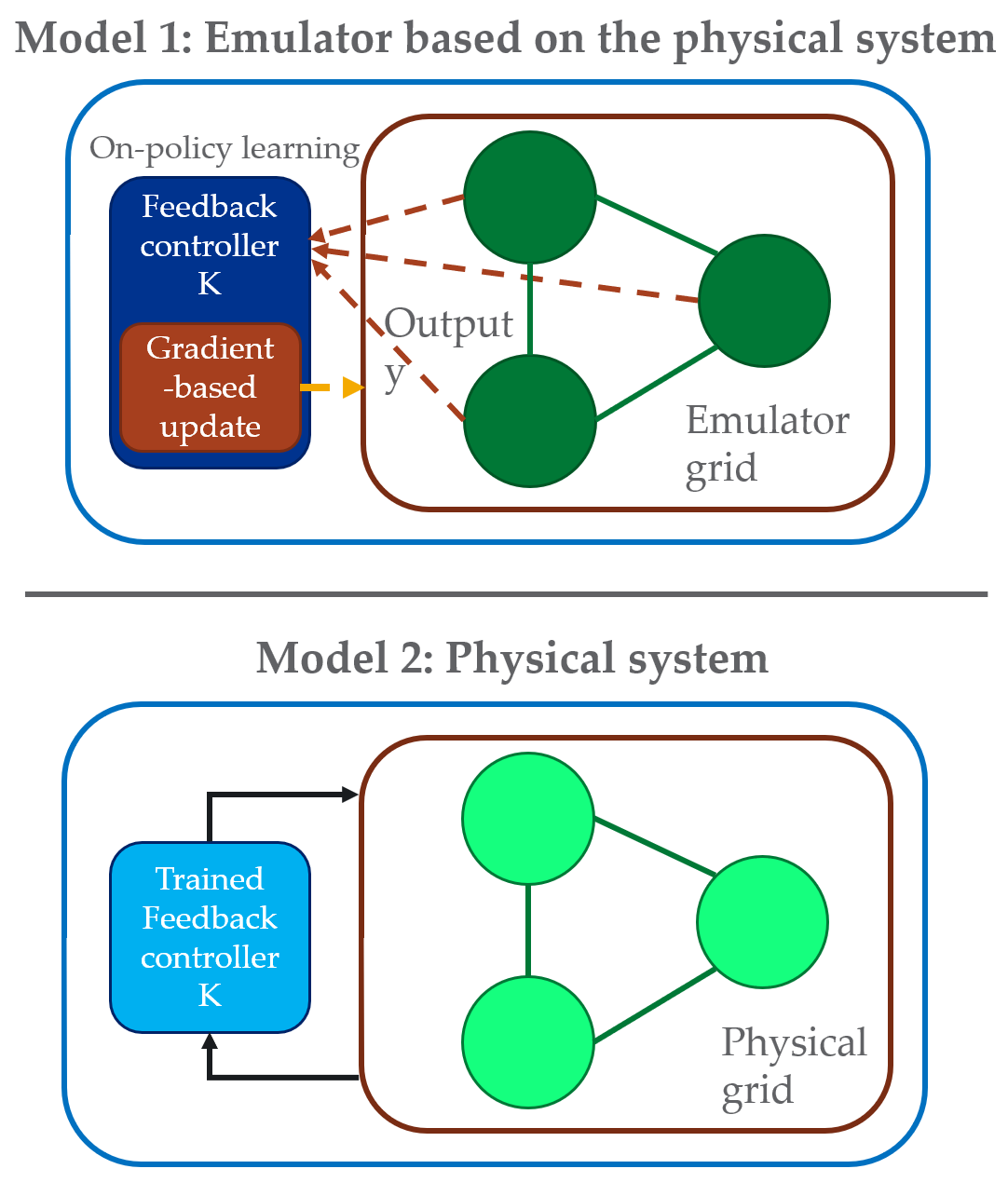}
	\caption{Problem setup with emulator (Model 1) and physical system (Model 2)}
	\label{fig:setup}
\end{figure}

The objective is to find an optimal feedback controller $K$ of the physical system, represented as Model~2 in Fig.~\ref{fig:setup}. {\KBB We also estimate the parameters of the system from the historical data and construct Model~1, which is the offline emulator based on the actual physical system.} Here, we formulate offline learning based on uncertainty data from the physical system. From the initial feedback gain $K^0$, we will observe the output $y$ with on-policy learning, meaning that the controller will be directly applied to the system with the uncertainty. The output result will be used to update $K$ by gradient-based update. By updating $K$ iteratively, the solution will approach the optimal feedback gain.

Note that the setup can be extended to the case with the nonlinear dynamics-based physical system in Model~2. Even though the system has nonlinear dynamics, we can construct Model~1 as a linear system by estimating the linear parameters and finding $K$ that can show the best result. The obtained $K$ may not be an optimal feedback gain due to the approximation from the nonlinear to linear dynamics, but it will show a decent result with much less computation burden by directly finding the feedback controller from the nonlinear formulation.

After we obtain the trained $K$ from Model~1, we adopt it to Model~2 to run on the physical system. Here, the robustness of the feedback gain is a key factor as the parameters between Model~1 and Model~2 may have errors due to the approximation. If the robustness is not high enough, the feedback controller will not work properly in Model~2 and thus needs a very precise approximation in Model~1. On the other hand, high robustness refers to stability in Model~2 even though there are some errors in approximations.
To formulate the algorithm, we start with the assumption of the emulator system with linear dynamics as:
\begin{align}
    &\dot{x} = Ax + Bu + w, x(0) = x_0\\
    & y = Cx
\end{align}
where $x\in \mathbb{R}^n, y \in \mathbb{R}^p, u \in \mathbb{R}^m$. We want to find static output feedback control:
\begin{align}
    u = -Ky = -KCx
\end{align}
that minimizes
\begin{align}
	R_0(K) &= \lim_{T \to \infty}\frac{1}{T}\mathbb{E} \sum_{t=0}^{T-1} [y^\top_t Q y_t + u^\top_t R u_t]  \label{eq:opt} 
\end{align}
considering the mean-variance risk constraint:
\begin{align}
	R_c(K) = \lim_{T \to \infty}\!\frac{1}{T} \mathbb{E} \sum_{t=0}^{T-1} \left(y_t^\top Q y_t - \mathbb{E}[y_t^\top Q y_t \vert h_t]\right)^2 \leq \delta  \nonumber 
\end{align}
Here, we consider the structured $K$ which is within the set $\mathcal{K}$ defined as
\begin{align}
	&\mathcal{K} = \{K: K_{a,b} = 0 \;\text{if and only if}\; (a,b) \notin \mathcal{E}) \}
\end{align}
where the structure pattern $\mathcal{E}$ is specified by the edges of a given communication or information-exchange graph. In other words, $K$ will have a sparse structure with the element as zero if the corresponding edges cannot exchange information.

By following \cite{tsiamis}, we can reformulate the constraint as
\begin{align}
	R_c(K) \!=\!  \lim_{T \rightarrow \!\infty} \frac{1}{T} \mathbb{E} \sum_{t=0}^{T-1} \left(4y_t^\top QCWC^\intercal Q^\intercal y_t+4y_t^\top Q M_3^y \right)\!\leq\!\bar{\delta}\! \label{eq:objr2}
\end{align}
with $\bar{\delta} = \delta-m_4^y + 4 \text{tr}\{(WQ_c)^2 \}$ and the (weighted) noise statistics given as
\begin{align}
	\bar{w}&=\mathbb{E}[w_t],\\
	W&=\mathbb{E}[(w_t-\bar{w})(w_t-\bar{w})^\top],\\ 
	M_3^y &= \mathbb{E}[(w_t-\bar{w})(w_t-\bar{w})^\top Q_c (w_t-\bar{w})],\\
	m_4^y&=\mathbb{E}[(w_t-\bar{w})^\top Q_c (w_t-\bar{w}) - \text{tr}(WQ_c)]^2.
\end{align} 
with $Q_c = C^\intercal Q C$. Note that \eqref{eq:objr2} no longer depends on the past trajectory $h_t$, and has a quadratic form which is the same form as the objective function \eqref{eq:opt}. 

To solve the optimization problem with the objective function as \eqref{eq:opt} and the risk-constraint as \eqref{eq:objr2}, we will formulate the Lagrangian function as below:
\begin{align}
\mathcal{L}(K, \lambda) =&R_0(K)+ \lambda [R_c(K)-\delta]\nonumber\\
 =&\lim_{T \rightarrow \infty} \frac{1}{T}\mathbb{E}\sum_{t=0}^{T-1} [y^\top_t (Q+4\lambda QCWC^\top Q^\top) y_t  \nonumber\\
 &+ 4\lambda y_t^\top Q M_3^y + u^\top_t R u_t]-\lambda\overline{\delta} \label{eq:lag}
\end{align}
By assuming that the problem is feasible and thus the dual variable $\lambda$ is finite, we can formulate the dual function with the bounded set $\mathcal{Y}:= [0,\Lambda]$ as
\begin{align}
\max_{\lambda \in \mathcal{Y}} \mathcal{D}(\lambda) = \max_{\lambda \in \mathcal{Y}}  \min_{K \in \mathcal{K}} \mathcal{L} (K, \lambda) \label{eq:dual}
\end{align}
Here, instead of \eqref{eq:dual}, we will consider the minimax counter problem:
\begin{align}
\min_{K\in\mathcal{K}}\Phi(K)~~\textrm{where} ~~\Phi(K):= \max_{\lambda \in \mathcal{Y}} \mathcal{L}(K, \lambda) \label{eq:dual_c}
\end{align}
Though \eqref{eq:dual} and \eqref{eq:dual_c} are different optimization problems, solving \eqref{eq:dual_c} will give the stationary point (SP) of \eqref{eq:dual} since the KKT conditions for \eqref{eq:dual} are related to the SP of \eqref{eq:dual_c}. Note that we can easily find $\lambda$ that maximizes $\mathcal{L}(K,\lambda)$ by choosing $\lambda = 0$ if the constraint is satisfied and $\lambda = \Lambda$ otherwise, since $\mathcal{L}(K,\lambda)$ is a linear function of $\lambda$.\\

\section{SGDmax Algorithm with ZOPG} \label{sec:ALG}

To solve \eqref{eq:dual_c}, we will adopt the SGDmax algorithm with ZOPG method to estimate the gradient \cite{KB}. Algorithm~\ref{alg:ZOO} shows how to find the zero-order policy gradient to update $K$ \cite{ZOO}. Specifically, we will choose the perturbation $U \in \mathcal{S}_\mathcal{K}$, with $\mathcal{S}_\mathcal{K}$ is the set of matrices that has the same structure in $\mathcal{K}$ with its magnitude being 1. After adopting the perturbation $rU$ with smoothing radius $r$ to current feedback gain $K$, the optimal $\lambda$ that maximizes the Lagrangian function value is chosen. By adopting $\lambda$, we can compute the estimated zero-order gradient. Note that the gradient is computed as $U$ multiplied by a scalar value, meaning that the estimated gradient follows the structure $\mathcal{K}$.
\begin{algorithm}[t]
	\SetAlgoLined
	\caption{Zero-Order Policy Gradient (ZOPG)}
	\label{alg:ZOO}
	\DontPrintSemicolon
	{\bf Inputs:} smoothing radius $r$, 
	the policy $K$ and its perturbation  $U \in \mathcal{S}_{\mathcal{K}}$, both of $n_\mathcal{K}$ non-zeros.\; 
	Obtain $\lambda' \leftarrow \argmax_{\lambda \in \mathcal{Y}} \mathcal{L}(K+rU, \lambda)$;\; 
Estimate the gradient $\hat{\nabla}_{\mathcal{K}} \mathcal{L} (K; U) = \frac{n_\mathcal{K}}{r} \mathcal{L}(K+rU, \lambda') U$.\;
	{\bf Return:} $\hat{\nabla}_{\mathcal{K}} \mathcal{L} (K; U)$.
\end{algorithm}
With the gradient from Algorithm~\ref{alg:ZOO}, we adopt the SGDmax algorithm introduced in Algorithm~\ref{alg:SGD} to find SP $K$ using gradient method. From the zero-order gradient for each $U_s$, the gradient is averaged and used to update $K^j$ in the iteration $j$ with gradient descent. Note that at the initialization, the set $\mathcal{G}^1$ and regarding parameters $L_0, \ell_0$ and $\rho_0$ are defined as
\begin{align}
&\mathcal{G}^1 := \{K\in\mathcal{K} |\Phi_{\mu_0}(K)\leq 10~\Phi_{\mu_0}(K^0)\}, \nonumber \\
&L_0\!:=\!\sup\limits_{K\in\mathcal{G}^0}{L_K}, \ell_0\!:=\!\sup\limits_{K\in\mathcal{G}^0}\ell_K. ~\textrm{and}~ \rho_0\!:=\!\inf\limits_{K \in \mathcal{G}^0} \min\{\beta_K, \psi_K\} \label{eq:const}
\end{align}
with Moreau envelope $\Phi_{\mu_0}(k)$ of $\Phi(K)$ defined by solving the minimization problem 
\begin{align}
    \Phi_{\mu_0}(K) := \min_{K'\in\mathcal{K}} \Phi(K') + \frac{1}{2{\mu_0}}\Vert K'-K\Vert^2, ~~\forall K \in \mathcal{K}, \label{eq:mor} 
\end{align}
where $\mu_0 = 1/(2\ell_0)$. Here, the sublevel set $\mathcal{G}^1$ ensures the Lipschitz and smoothness properties for any $K \in \mathcal{G}^1$.{\KB Note that smoothing radius $r$, stepsize $\eta$, and the minimum number of iteration $J$ are determined by following \cite[Theorem 2]{KB} which is refered here as Lemma~\ref{lem:proof}.

\begin{lemma}
With an initial $K^0 \in \mathcal{K}$ and a given $\epsilon > 0$, by setting the parameters as 
	\begin{align}
		&r \leq \min\Big\{\rho_0, \frac{L_0 \sqrt{M}}{\ell_0} \Big\}, \; \eta \leq \frac{\epsilon^2}{\alpha\ell_0(L_0^2+\ell_0^2r^2/M)},\nonumber\\
		\textrm{and}~&J = \frac{2\sqrt{10\alpha}\Phi_{\mu_0}(K^0)}{\eta \epsilon^2}\label{eq:thm2}
	\end{align}
with $\mathcal{G}^1, L_0, \ell_0$ and $\rho_0$ following \eqref{eq:const}, Algorithm~\ref{alg:SGD} converges to $\epsilon$-SP $K_\epsilon$ that satisfies $\Vert \Phi_{\mu_0}(K_\epsilon)\Vert \leq \epsilon$, with probability of at least $(0.9-\frac{4}{\alpha}-\frac{4}{\sqrt{10\alpha}})$. 
\label{lem:proof}
\end{lemma}
}
\begin{algorithm}[t]
	\SetAlgoLined
	\caption{Stochastic Gradient Descent with max-oracle (SGDmax)}
	\label{alg:SGD}
	\DontPrintSemicolon
	{\bf Inputs:} A feasible and stable policy $K^0$, upper bound $\Lambda$ for $\lambda$, threshold $\epsilon$, and number of ZOPG samples $M$.\; 
		Determine $L_0, \ell_0$, and $\rho_0$ with the set $\mathcal{G}^1$ and compute $r$, $\eta$, and $J$; \;
	\For{$j = 0, 1, \ldots, J-1$}{
		\For{$s=1,\ldots,M$}{
		Sample the random $U_s \in \mathcal{S}_\mathcal{K}$; \;
		Use Algorithm~\ref{alg:ZOO} to return $\hat{\nabla}\mathcal{L}_\mathcal{K}(K^j;U_s)$. 
		}
		Update $K^{j+1} \leftarrow K^{j}-\eta \left(\frac{1}{M} \sum_{s=1}^M \hat{\nabla} \mathcal{L}(K^j;U_s)\right)$.}
	{\bf Return:} the final iterate $K^{J}$.
\end{algorithm}

{\KB \noindent Proof Sketch: Similar to \cite{KB}, the sketch of the proof for Lemma~\ref{lem:proof} is as follows.} From the result of recent works \cite{bu, malik, harvard}, the Lipschitz and smoothness properties of $\mathcal{L}(K,\lambda)$ are still valid by \cite[Lemma~1]{KB}. {\KBB Note that the output feedback involves only the linear combination from the original states, thereby such functional properties are maintained.} Therefore, we can bound the expectation of $\Phi_{\mu_0}(K^{j+1})$ from the expectation of $\Phi_{\mu_0}(K^{j})$ as below;
\begin{align}
	\mathbb{E}&\left[\Phi_{\mu_0}(K^{j+1})\right]\leq \mathbb{E}\left[\Phi_{\mu_0}(K^{j})\right]\nonumber\\
	& - \frac{\eta}{4} \mathbb{E} \Vert \nabla \Phi_{\mu_0}(K^{j})\Vert^2 + \eta^2 \ell_0 \left(L_0^2+{\ell_0^2 r^2}/{M}\right) \label{eq:pbd2}
\end{align}
Here, we can claim the convergence to an SP by showing these two aspects. First, we will verify the SP condition, i.e. $\Vert \nabla \Phi_{\mu_0}(K) \Vert \leq \epsilon$ in the expectation sense as the iteration continues. Second, we will verify that the probability that the algorithm fails to converge to the SP is bounded.
First, we can verify the SP condition by summing up \eqref{eq:pbd2} over iterations $j=0, ..., J-1$,
\begin{align}
    &\frac{1}{J}\sum_{j=0}^{J-1}\mathbb{E}\left[\Vert \nabla\Phi_{\mu_0}(K^j)\Vert^2\right] \nonumber\\
    \leq&\frac{4\left[\Phi_{\mu_0}(K^0)\!-\!\mathbb{E}[\Phi_{\mu_0}(K^J)]\right]}{J\eta}+\!4\eta\ell_0(L_0^2+\ell_0^2r^2/M),
\end{align}
which gives the upper bound on order of $\epsilon^2$ by choosing the stepsize $\eta=O(\epsilon^2)$. 
Second, the probability that $\{K^j\}$ exceeds $\mathcal{G}^1$ is represented as in \eqref{eq:prob},
\begin{align}
    &\mathbb{P}\left( \frac{1}{J} \sum_{j=0}^{J-1} \Vert \nabla\Phi_{\mu_0}(K^j)\Vert^2 \geq \epsilon^2 \right) \nonumber\\
    = & \mathbb{P}\left( \frac{1}{J} \sum_{j=0}^{J-1}\Vert \nabla\Phi_{\mu_0}(K^j)\Vert^2 \geq \epsilon^2, \tau > J \right) \nonumber\\
    + &\mathbb{P}\left( \frac{1}{J} \sum_{j=0}^{J-1}\Vert \nabla\Phi_{\mu_0}(K^j)\Vert^2 \geq \epsilon^2, \tau \leq J \right), \label{eq:prob}
\end{align}
where $\tau := \min\{j\geq 0: K^j \notin \mathcal{G}^1 \}$, meaning the iteration index that $\{K^j\}$ exceeds $\mathcal{G}^1$.
Here, we can bound the first and second terms on the right-hand side by following the Markov's inequality and the Doob's maximal inequality with the supermartingale property. The detailed proof can be found at \cite[Appendix~B]{KB}.

\section{Simulation Study} \label{sec:SR}

\begin{figure}[t]
	\centering
	\includegraphics[width=0.65\linewidth]{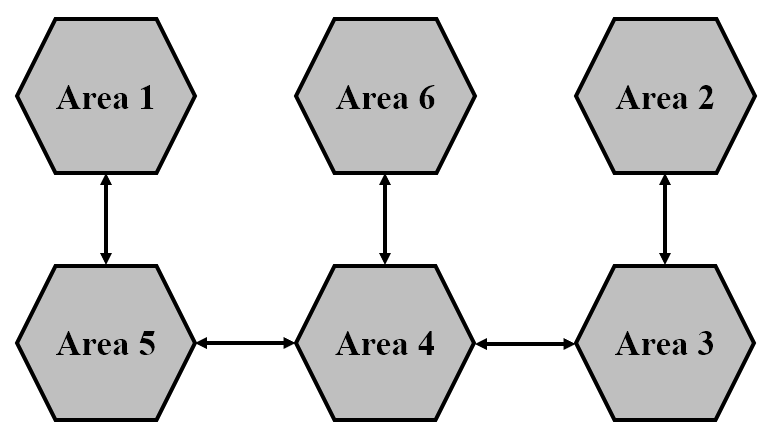}
	\caption{Power grid with six interconnected areas.}
	\label{fig:network}
\end{figure}

To show the effectiveness of the proposed method, we consider the LFC in the multi-area power system. Assume there are $N$ number of radially connected areas through tie lines depicted in Fig.~\ref{fig:network} in case $N=6$ as an example. The corresponding parameters are indicated in \cite[Table I]{DLQR}.
Assume the information-exchange graph $\mathcal{G}=(\mathcal{V},\mathcal{E})$ follows the physical interconnection between areas in Fig.~\ref{fig:network}. Each node $i\in \mathcal{V}$ represents an area and an edge $(i,j) \in \mathcal{E}$ between two areas indicates that they are connected through a tie-line. 

\begin{figure}[t]
	\centering
	\includegraphics[width=\linewidth]{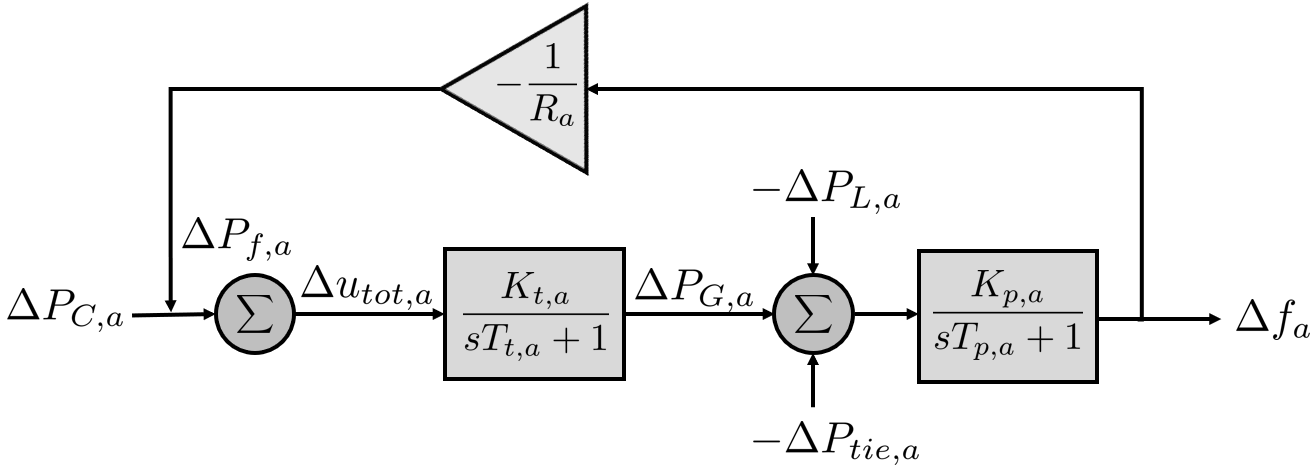}
	\caption{Block representation of load frequency control in $i$-th area \cite{KB}}
	\label{fig:dynamics}
\end{figure}

For the aggregated generators in the area, the open-loop linearized dynamics of the $i$-th area are depicted in Fig.~\ref{fig:dynamics} \cite{bevrani, KB}.
Notation $\Delta$ implies the deviation from the steady-state operation. When there is a frequency deviation $\Delta f_i$ happens in $i$-th area, the primary frequency control action $\Delta P_{f,i}$ is made proportional to $\Delta f_i$, as $\Delta P_{f,i} = -(1/R_i) \Delta f_i$. As a secondary control, automatic generation control (AGC) signal $\Delta P_{C,i}$ is determined based on the control design. The total control signal $\Delta u_{tot, i}$ combines primary and secondary control of generators, and eventually construct the signal total power output of generators denoted as $\Delta P_{G,i}$. Third, the unknown power load deviation $\Delta P_{L,i}$ and the total power inflow $\Delta P_{tie,i}$ are added to $\Delta P_{G,i}$ and affects the frequency deviation $\Delta f_i$ of the next time slot. Note that $\Delta P_{tie,i}$ is the summation of power inflow from the interconnected areas that are proportional to the difference in frequency as
\begin{align}
	\Delta P_{tie,i} = \sum_{j\in \mathcal{N}_i} K_{tie,i,j} (\Delta f_i - \Delta f_j)\Delta t,
\end{align}   
{\KB where $j\in\mathcal{N}_i$ indicates the set of areas that are connected to $i$-th area, $\Delta t$ means the discrete time interval, and $K_{tie}$ refers to a synchronization coefficient between area $i$ and $j$.} 
In addition to the dynamics indicated in Fig.~\ref{fig:dynamics}, we adopt the integral action of Area Control Error (ACE) $z_i$ defined as $z_i = \beta_i \Delta f_i + \Delta P_{tie,i}$, with bias factor $\beta_i$ often being defined as $D_i + 1/R_i$ \cite{DLQR_en}.   

The state, control variable, and the unknown random disturbance of $i$-th area as $x_i = [\Delta f_i, \Delta P_{G,i}, \Delta P_{tie,i}, \int z_i]^\top$, $u_t = \Delta P_{C,i}$ and $w_i = \Delta P_{L,i}$ respectively. With an assumption that all areas are identical, the aggregated dynamics of the network can be represented as linear dynamics after removing the subscript $i$:
\begin{align}
	\dot{\tilde{x}} = (I_N \otimes A_1 + \mathcal{L} \otimes A_2)\tilde{x} + (I_N \otimes B_u) \tilde{u} + (I_N \otimes B_w)\tilde{w} \label{eq:agg_ss}
\end{align}

\noindent where $\tilde{x}\!=\!\text{col}(x_1,\cdots,x_N), \tilde{u}\!=\!\text{col}(\Delta P_{C,1},\cdots\!,\Delta P_{C,N}), \tilde{w} = \text{col}(\Delta P_{L,1},\cdots,\Delta P_{L,N})$ and the matrices $A_1, A_2, B_u$ and $B_w$ follow \cite[Section~V]{KB}. 


Lastly, the output feedback $\tilde{y}$ is computed as $\tilde{y}=C \tilde{x}$, with $C$ indicating the output that can be observed. Here, since the frequency deviation $\Delta f_i$ is hard to track when considering the aggregated generators in an area, we will observe only the sum of total power output of generators and total power inflow, i.e. $\Delta P_{G,i}$ and $\Delta P_{tie,i}$. Later, we show an experiment when we include the area frequencies in the feedback. Also, we experiment with the design without using the conventional ACE inputs, however, indirectly, the tie-line power flow captures their impact. \color{black} Accordingly, $C$ represents the output that can be observed based on the information-exchange graph. {\KB In Fig.~\ref{fig:network}, the corresponding matrix $C$ will be a $6 \times 24$ matrix, with the elements $C_{i,4j-2}$ and $C_{i,4j-1}$ being 1 if the areas $i$ and $j$ are connected. 
Here, each column represents the provided information from the neighborhood areas. For example, area~2 can use the information from area~3 and itself, which corresponds to column 2 that has $C_{2,6}, C_{2,7}, C_{2,9}$ and $C_{2,10}$ as 1.} Here, we adopt the linear control policy $\tilde{u} = -K \tilde{y} = -KC\tilde{x}$, with $K$ following the sparse structure $\mathcal{K}$. It implies that $K_{i,j} = 0$ when node $i$ and $j$ are not connected in the information-exchange graph.

\subsection{Training Result with Model~1}
\begin{figure}[t]
	\centering
	\includegraphics[width=0.95\linewidth]{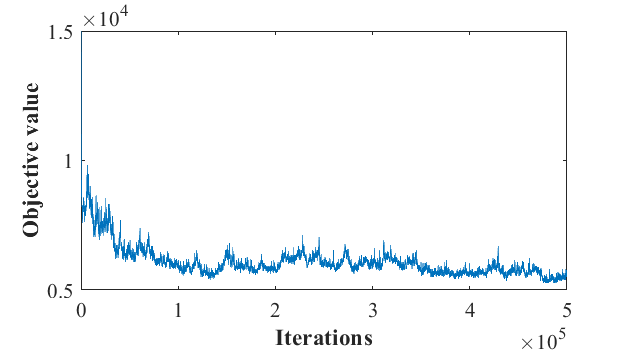}
	\caption{LQR objective trajectory.}
	\label{fig:obj}
\end{figure}

\begin{figure}[t]
	\centering
	\includegraphics[width=0.95\linewidth]{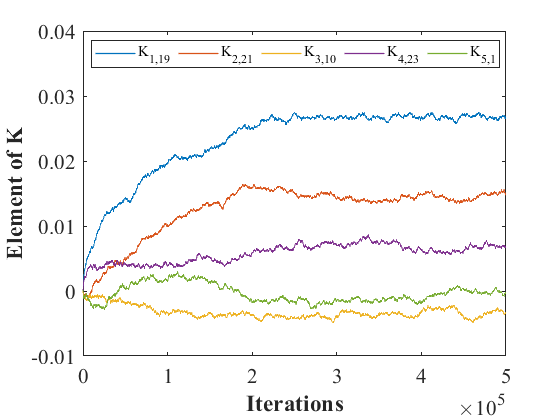}
	\caption{Updates of selected elements of $K$ during learning.}
	\label{fig:element}
\end{figure}

To find the optimal feedback gain $K$, Algorithm~\ref{alg:SGD} is adopted with a stepsize $\eta = 10^{-4}$, smoothing radius $r=0.1$, and $M=100$ samples for ZOPG. Figure~\ref{fig:obj} shows that LQR objective value converges to a steady-state with a sufficient number of updates. 
In addition, Fig~\ref{fig:element} depicts the trajectory of selected elements of $K$ during the updates. Even though there is some fluctuation due to the estimated gradient, all values converge as iteration continues, meaning that $K$ reaches the SP of the objective function.
\subsection{Testing Result with Model~2}
\begin{figure}[t]
	\centering
	\includegraphics[width=0.95\linewidth]{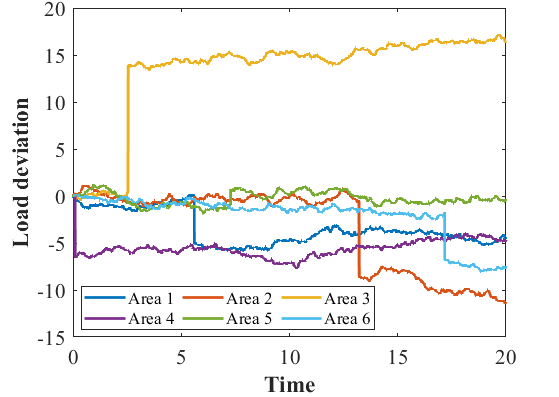}
	\caption{Load deviations of areas}
	\label{fig:load}
\end{figure}

\begin{figure}[t]
	\centering
	\includegraphics[width=0.95\linewidth]{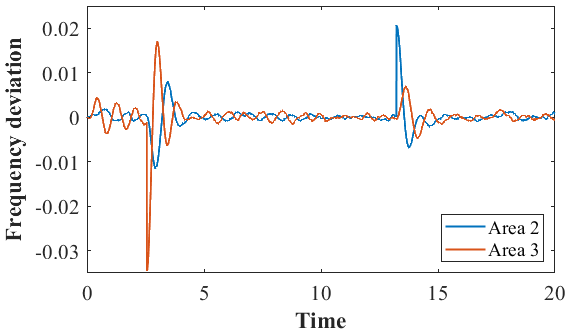}
	\caption{Frequency deviation of Area~2 and Area~3}
	\label{fig:freq}
\end{figure}

\begin{figure}[t]
	\centering
	\includegraphics[width=0.95\linewidth]{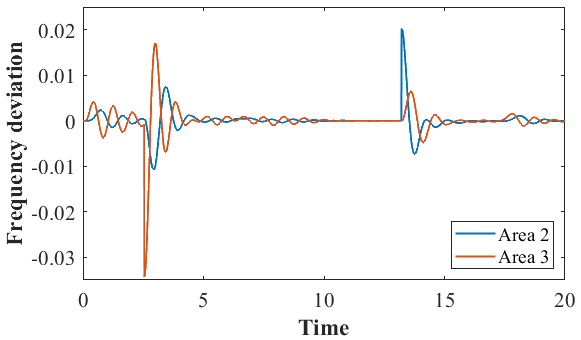}
	\caption{Frequency deviation of Area~2 and Area~3 when including frequency deviation in the feedback}
	\label{fig:freq2}
\end{figure}

With the feedback gain $K$ obtained from the training, test has been made by creating a scenario in a 20-second window. Fig~\ref{fig:load} shows one of the scenario on the changes in the load in each area. To show the effectiveness of the control, each area has one large random changes during a 20-second time period. For example, area~3 has huge increase in the load around 3 seconds according to Fig~\ref{fig:load}. 
{\KB Fig~\ref{fig:freq} 
shows the frequency deviation 
of Area~2 and Area~3 with the obtained $K$ from the training. Note that Area~2 and Area~3 are interconnected according to Fig.~\ref{fig:network} with output values representing the sum of total power output of generator and total power inflow. As seen in Fig~\ref{fig:freq}, the frequency changes dramatically when there is a huge change in load as indicated in Fig.~\ref{fig:load}. Frequency deviation is stabled to zero in a short time period meaning that the feedback gain $K$ can successfully control the frequency from the uncertain load changes. Here, the connected area is also affected since they are connected to each other and thus affects the control decision based on the structure of $K$.} {\KBB Moreover, the performance in Fig~\ref{fig:freq} is similar to  Fig.~\ref{fig:freq2}, which includes the frequency deviation in the feedback. It shows that our proposed output feedback works quite well even when we do not have the information on frequency, with the influence of $\Delta P_{G,i}$ and $\Delta P_{tie,i}$ captured through peer-to-peer information exchange.}

\subsection{Robustness Test Result}
As discussed, the parameters of the dynamics in Model~1 may have errors compared to Model~2 due to the estimation or the approximation from the nonlinear dynamics to linear dynamics. Therefore, checking out the robustness of the obtained feedback gain is crucial when adopting the obtained feedback gain to the physical system.

\begin{figure}[t]
	\centering
	\includegraphics[width=0.95\linewidth]{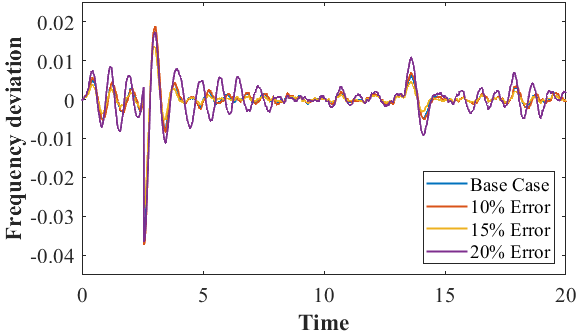}
	\caption{Frequency deviation of Area~3 for robustness test with Model 1}
	\label{fig:robust}
\end{figure}

To test the robustness, we observe the frequency deviation of Area~3 with the changes in the parameters by 10\%, 15\% and 20\% of their actual values, respectively. Fig~\ref{fig:robust} shows the result of four cases. As seen in the figure, the controller successfully controls the frequency deviation until 15\% error, but the fluctuation increases in 20\% error case and thus takes much more time to reach the stable frequency. The result implies that even when the estimated parameters in Model~1 have errors with Model~2, the feedback gain from Model~1 can still perform a good frequency control in the physical system in around 15\% of the error. 

\section{Conclusions} \label{sec:CON}
The paper presented a learning-based method to solve the P2P-based LFC problem, while considering three challenges including model uncertainty, the mean-variance risk constraint in load uncertainty, and output structured feedback. To train a feedback controller based on the data from a physical grid, an emulator grid that approximates the dynamics of the physical grid was formulated and the learning is performed in the emulator grid. Here, starting from the constrained LQR problem with output structured feedback, we reformulated the risk constraint to have the same form as the LQR objective function and construct the dual problem as a minimax problem. By adopting SGDmax algorithm with an estimated ZOPG, we obtained the trained output feedback that we can adapt to a physical grid with a robustness test. Numerical tests with radially interconnected areas in a power system showed that the feedback controller from learning can successfully control the frequency of the areas while satisfying the risk constraints. In addition, robustness test indicated that in the current setting of the model the obtained feedback controller can be adapted to a physical system with 20\% robustness, meaning that the controller obtained from the approximated emulator grid can be directly adopted to the physical system. {\KBB Future work includes uncertainty from distributed energy resources such as wind/solar or EVs, and adopts the proposed method to real power system models.}

\bibliography{bibliography}
\bibliographystyle{IEEEtran}
\itemsep2pt

\end{document}